\documentclass [12pt]{article}
\newfont{\bbb} {msbm10}
\newcommand{\Bbb}[1]{\mbox{\bbb#1}}
\newcommand{\bS}{\Bbb{S}}

\newcommand{\Z}{\Bbb{Z}}
\newcommand{\Q}{\Bbb{Q}}

\newcommand{\sm}{\setminus}
\newcommand{\sbs}{\subset}
\newcommand{\ra}{\rightarrow}

\newcommand{\p}{\partial}

\newcommand{\cB}{{\cal{B}}}

\newcommand{\uc}{\underline{c}}
\newcommand{\ua}{\underline{a}}
\newcommand{\ub}{\underline{b}}

\newcommand{\bG}{{\bf G}}
\usepackage[all]{xy}

\begin{document}

\title{Branched Covers of Hyperbolic Manifolds and Harmonic Maps}
\author{F. T. Farrell and P. Ontaneda\thanks{The first author was
partially supported by a NSF grant.
The second author was supported in part 
by a research  grant from CAPES, Brazil.}}
\date{}

\maketitle

\noindent {\bf \large 1. Introduction and statements of the results}
\vspace{.2in}

Let $f:M\ra N$ be a homotopy equivalence between closed negatively curved manifolds.
The fundamental existence results of Eells and Sampson \cite{ES} and uniqueness of
Hartmann \cite{Har} and Al'ber \cite{A'} grant the existence of a unique harmonic
map $h$ homotopic to $f$.
Based on the enormous success of the harmonic map technique Lawson and Yau conjectured that
the harmonic map $h$ should be a diffeomorphism. This conjecture was proved to be false by 
Farrell and Jones \cite{FJ1} in every dimension in which exotic spheres exist.
They constructed examples of homeomorphisms  $f:M\ra N$ between closed negatively curved manifolds
for which $f$ is not homotopic to a diffeomorphism.
These counterexamples were later obtained also in dimension six by Ontaneda \cite{O} and
later generalized by Farrell, Jones and Ontaneda to all dimensions $>5$ \cite{FJO1}. In fact, in \cite{O} and
\cite{FJO1} examples are given for which $f$ is not even homotopic to a
$PL$ homeomorphism. The fact that $f$ is not homotopic to a $PL$ homeomorphism has several 
interesting strong consequences that imply certain limitations of well known powerful analytic
methods in geometry \cite{FJO2}, \cite{FOR}, \cite{FO1}, \cite{FO2} (see \cite{FO3} for a survey).\\

In all the examples mentioned above one of the manifolds is always
a hyperbolic manifold. Hence, both manifolds $M$ and $N$ have the homotopy type of a hyperbolic
manifold (hence the homotopy type of a locally symmetric space). We call these 
{\it examples of the first kind}.\\

In \cite{Ar} Ardanza also gave counterexamples to the Lawson-Yau conjecture.
In his examples the manifolds $M$ and $N$ are not homotopy equivalent to a locally symmetric space;
in particular, they are not homotopy equivalent to a hyperbolic manifold. We call these {\it examples of the
second kind}.
In these examples the map $f$ is not homotopic to a diffeomorphism and exist in dimensions
$4n-1$, $n\geq 2$. Hence Ardanza's result is an analogue of Farrell-Jones result \cite{FJ1} for
examples of the second kind.
His constructions use branched covers of hyperbolic manifolds. Recall that 
Gromov and Thurston \cite {GT} proved that large branched covers of hyperbolic manifolds do not
have the homotopy type of a locally symmetric space. \\

In this paper we show that most of the results obtained for examples of the first kind 
in \cite{O}, \cite{FJO1}, \cite{FJO2}, \cite{FOR}, \cite{FO1}, \cite{FO2}, \cite{FO3} are also
true for examples of the second kind. We now state our main results.\\

First we extend Ardanza's result to every
dimension $>5$ and also with the stronger property that $f$ is not even homotopic to a $PL$
homeomorphism. This is an analogue of the result in \cite{FJO1} for examples of the second kind.\\

\noindent {\bf Theorem 1.} {\it Given $\epsilon >0$ and  $n>5$ there
are closed Riemannian manifolds $M_1$, $M_2$ and a homotopy equivalence
$f:M_1\ra M_2$ such that:} 
\begin{enumerate}
\item[{1.}] {\it $f$ is not homotopic to a $PL$ homeomorphism.
In particular, the unique harmonic map homotopic to $f$ is not a diffeomorphism.}

\item[{2.}] {\it $M_1$ and $M_2$ have $\epsilon$-pinched to $-1$ sectional curvatures.}

\item[{3.}] {\it $M_1$ and $M_2$ do not have the homotopy type of a locally symmetric space.
In particular, they are not homotopy equivalent to a hyperbolic manifold.}
\end{enumerate}
\vspace{.2in}

\noindent {\bf Remark.} This result is a little weaker than the result in \cite{Ar} 
(for examples of the second kind) and than the one in \cite{FJO1} (for examples of the first kind):
in \cite{Ar}  $M_1$ is not diffeomorphic to $M_2$ and in \cite{FJO1} $M_1$ is not $PL$ homeomorphic to $M_2$.
The Theorem above states the existence of a particular map $f$ that is not homotopic to a 
$PL$ homeomorphism but we do not know whether $M_1$ and $M_2$ 
are $PL$ homeomorphic. The missing ingredient is ``differentiable rigidity" which, for  
examples of the first kind, is provided by Mostow's Rigidity Theorem:
every homotopy equivalence between hyperbolic manifolds is homotopic to a diffeomorphism
(in fact an isometry). We do not know whether every homotopy equivalence between branched
covers of hyperbolic manifolds is homotopic to a diffeomorphism.\\

As for the case of examples of the first kind,
the fact that the map $f$ is not homotopic to a $PL$ homeomorphism has the following 
interesting consequence:\\

\noindent{\bf Theorem 2.}  {\it Given $\epsilon >0$ and  $n>5$ there are
closed Riemannian manifolds $M_1$, $M_2$ and a harmonic homotopy equivalence
$h:M_1\ra M_2$ such that:} 
\begin{enumerate}
\item[{1.}] {\it $h$ is not one-to-one.}

\item[{2.}] {\it $M_1$ and $M_2$ have $\epsilon$-pinched to $-1$ sectional curvatures.}

\item[{3.}] {\it $M_1$ and $M_2$ do not have the homotopy type of a locally symmetric space.
In particular, they are not homotopy equivalent to a hyperbolic manifold.}
\end{enumerate}
\vspace{.2in}

This Theorem can be directly deduced from Theorem 1 and the $C^\infty -$ Hauptvermutung of Scharlemann
and Siebenmann \cite{SC1}.
Also, if Poincare's conjecture for three dimensional manifolds is true then the map $h$ above is not even a 
cellular map, hence it cannot be approximated by homeomorphisms. The proof of this
fact is similar to the proof for examples of the first kind. For more details see
\cite{FJO2} or \cite{FO3}.\\

In the examples provided by the Theorem above, the main obstruction to $h$ being a diffeomorphism 
or a homeomorphism is that $h$ is not homotopic to a $PL$ homeomorphism.  
We may ask then what happens if this
obstruction vanishes, that is, if $h$ is homotopic to a $PL$ homeomorphism or even 
homotopic to a diffeomorphism.
This was considered in  Problem 111 of the list compiled by S.-T.  
Yau in  \cite{Yau}.  Here is a restatement of this problem.\\

\noindent {\bf Problem 111 of \cite{Yau}}.  Let $f:M_1\rightarrow M_2$ be a
diffeomorphism between two compact manifolds with negative curvature.  If
$h:M_1\rightarrow M_2$ is the unique harmonic map which is homotopic to $f$, is $h$ a
homeomorphism?, or equivalently, is $h$ one-to-one?\\

(This problem had been reposed in \cite{Sta} as Grand Challenge Problem 3.6.)  The
answer to the problem was proved to be yes when dim$M_1=2$ by Schoen-Yau
\cite{SY} and Sampson  \cite{Sa}.  But it was proved by Farrell, Ontaneda and Raghunathan
\cite{FOR} that the answer to this question is in general negative for dimensions $>5$. The 
counterexamples constructed in \cite{FOR} are examples of the first kind. Here we construct 
counterexamples of the second kind:\\

\noindent {\bf Theorem 3.} {\it For every integer $n>5$ and $\epsilon >0$, there is a diffeomorphism $f:M_1
\rightarrow M_2$ between a pair of closed $n$-dimensional Riemannian manifolds such that:} 
\begin{enumerate}
\item[{1.}] {\it The unique harmonic map homotopic to $f$ is not a diffeomorphism.}

\item[{2.}] {\it $M$ and $N$ have $\epsilon$-pinched to $-1$ sectional curvatures.}

\item[{3.}] {\it $M$ and $N$ do not have the homotopy type of a locally symmetric space.
In particular, they are not homotopy equivalent to a hyperbolic manifold.}
\end{enumerate}
\vspace{.2in}

Also, in \cite{FO1} we constructed examples of harmonic maps $h:M\ra N$ between 
$\epsilon$-pinched to $-1$ closed  Riemannian manifolds such that
$h$ can be approximated by diffeomorphisms, but $h$ is not a diffeomorphism
(in particular $h$ is a cellular map).  The examples in \cite{FO1} are
examples of the first kind. Here we also construct counterexamples of the second
kind:\\

\noindent {\bf Theorem 4.} {\it For every integer $n>10$ and $\epsilon >0$, there is a harmonic map $h:M_1
\rightarrow M_2$ between a pair of closed $n$-dimensional Riemannian manifolds such that:} 
\begin{enumerate}
\item[{1.}] {\it  The harmonic map $h$ is not a diffeomorphism.}

\item[{2.}] {\it The harmonic map $h$ can be approximated by diffeomorphisms (in the $C^{\infty}$ topology).}

\item[{3.}] {\it $M_1$ and $M_2$ have $\epsilon$-pinched to $-1$ sectional curvatures.}

\item[{4.}] {\it $M_1$ and $M_2$ do not have the homotopy type of a locally symmetric space.
In particular, they are not homotopy equivalent to a hyperbolic manifold.}
\end{enumerate}
\vspace{.2in}

As for examples of the first kind, if Poincare's Conjecture for three dimensional manifolds is true
we can say a little more in this case (see \cite{FO1} or \cite{FO3}).
Also, in Theorems 1,2, and 3 above (and maybe also in Theorem 4) we can replace the word
``harmonic" by ``natural". The concept of natural map was defined by
G. Besson, G. Courtois and S. Gallot (e.g. see \cite{BCG}). These maps have very powerful
dynamic and geometric properties and are also used to prove rigidity results.
(For more details see \cite{FO3}).\\

In \cite{FO2} the results of \cite{FO1} were used to construct examples of
$\epsilon$-pinched to $-1$ closed  Riemannian manifolds for which the Ricci flow
does not converge smoothly. These examples are examples of the first kind.
Here we show also that the constructions used to prove Theorem 4 (which are analogous
to the constructions in \cite{FO1}) can also be used to produce 
examples of the second kind of  $\epsilon$-pinched to $-1$ closed  Riemannian manifolds 
for which the Ricci flow does not converge smoothly.\\

\noindent {\bf Remark.}
We say the the Ricci flow for a negatively
curved Riemannian metric $h$ {\it converges smoothly}  if the Ricci flow, starting at $h$, 
is defined for all $t$ and converges (in the $C^\infty$ topology) to a well defined 
negatively curved (Einstein) metric.
The next Theorem shows the existence of pinched negatively curved metrics for which the
Ricci flow does not converge smoothly. \\

\noindent {\bf Theorem 5.} {\it For every integer $n>10$ and $\epsilon >0$, there is a  
closed smooth $n$-dimensional Riemannian manifold $N$ such that:} 
\begin{enumerate}
\item[{1.}] {\it $N$ admits a Riemannian metric with sectional curvatures in
$[-1-\epsilon ,-1+\epsilon]$ for which the Ricci flow does not converge smoothly.}

\item[{2.}] {\it  $N$ does not have the homotopy type of a locally symmetric space.
In particular, $N$ is not homotopy equivalent to a hyperbolic manifold.}
\end{enumerate}
\vspace{.2in}

The proofs of all Theorems above use the following Proposition:\\

\noindent {\bf Proposition.} {\it For every integer $n>5$ and $r,s >0$, there 
are closed connected orientable Riemannian manifolds  $M$, $N$, $T$, $P$  such that:} 
\begin{enumerate}
\item[{1.}] {\it  $M$ is a  $n$-dimensional hyperbolic manifold,
$N$ and $P$ are totally geodesic $(n-1)$-dimensional submanifolds of $M$ and 
$T$ is a totally geodesic $(n-2)$-dimensional submanifold of $M$. $N$ intersects $T$ and $P$ transversally.}

\item[{2.}] {\it The isometry class of $N$ does not depend on $r$ (only on $n$ and $s$).}

\item[{3.}]  {\it $0\neq [N\cap T]\in H_{n-3}(M,\Z_2 )$, where  $[N\cap T]$  is the $\Z_2$-homology class 
represented by the $(n-3)$-dimensional submanifold $N\cap T$.}

\item[{4.}] {\it The normal geodesic tubular neighborhood of $N$ has width $>r$ and
the normal geodesic tubular neighborhood of $N\cap P$ has width $>s$.}
\end{enumerate}
\vspace{.2in}

{\bf Remark.} The Proposition remains valid if we replace $N$ in item 2 by $T$.
(We cannot choose {\it both} $N$ and $T$ to be independent of the widths of their
normal geodesic tubular neighborhoods.)
Also, we can have $T\sbs P$, if we choose so.\\

In the next section we prove the Proposition and show how it implies the Theorems.
The proofs of the Theorems resemble the proofs of the  corresponding 
Theorems for the case of examples of the first kind presented in 
\cite{O}, \cite{FJO1}, \cite{FJO2}, \cite{FOR}, \cite{FO1}, \cite{FO2}, \cite{FO3}.
We will refer to these papers.\\

We are grateful to J-F. Lafont for his help.\\
\vspace{.4in}

\noindent {\bf \large 2. Proofs of the results}
\vspace{.2in}

First we prove the Proposition.\\

\noindent {\bf Proof of the Proposition.} We use all notation from \cite{FOR}.
Let $\bG$, $\bG_1$, $\bG_2$ be as in the proof of lemma of \cite{FOR}, p. 233.
We have that $\bG_i=R_{k/\Q}\bG_i'$ where $\bG_i'$ denotes the $k$-algebraic
group $SO(f_i)$ and $f_i$ is the restriction of the form $f$ to the subspace
of $E$ generated by $\cB_i$, see\cite{FOR} p.245. To be specific, choose 
$\cB_1=\cB\sm\{ e_1\}$ and $\cB_2=\cB\sm\{ e_2, e_3\}$.\\

For ideals $\ua$, $\ub$, $\uc\,$ let $\Phi (\ua , \ub, \uc )$ be the arithmetic
subgroup of $\bG$ constructed in p. 234 of \cite{FOR}.
We will need the following two properties of $\Phi$. These properties
can easily be checked directly from the definition.
\begin{enumerate}
\item[{a.}] $\Gamma (\uc )\sbs\Phi (\ua , \ub, \uc )\sbs \Gamma (\ua )$

\item[{b.}] If $\ua '\sbs\ua$, $\ub '\sbs\ub$, $\uc '\sbs\uc$ then 
$\Phi (\ua ', \ub ', \uc ')$ is a subgroup of finite index of $\Phi (\ua , \ub, \uc )$.
\end{enumerate}

The following is shown in \cite{FOR}:\\

There is an ideal $\ua_0$ of $\Z$ such that for every ideal $\ua\sbs\ua_0$ there is an 
ideal $\ub (\ua )\sbs\ua$ with the following properties. For any ideal $\uc\sbs\ub (\ua )$
define $\Phi =\Phi ( \ua , \ub (\ua ), \uc )$, and $\Phi_i$, $M$, $N$, $T$, as in p. 234
of \cite{FOR}. Then $M$, $N$ and $T$  are closed connected orientable manifolds that satisfy:
\begin{enumerate}
\item[{1.}]  $M$ is a  $n$-dimensional hyperbolic manifold,
$N$ is a totally geodesic $(n-1)$-dimensional submanifold of $M$ and 
$T$ is a totally geodesic $(n-2)$-dimensional submanifold of $M$. $N$ and $T$ intersect transversally.

\item[{2.}] The isometry class of $N$ does not depend on $\uc$.

\item[{3.}]  $0\neq [N\cap T]\in H_{n-3}(M,\Z_2 )$, where  $[N\cap T]$  is the $\Z_2$-homology class 
represented by the $n-3$-dimensional submanifold $N\cap T$.
\end{enumerate}

Moreover, given $r$ we can choose $\uc$  such that 
\begin{enumerate}
\item[{4.}]  The normal geodesic tubular neighborhood of $N$  has width $>r$.
\end{enumerate}

We have to define $P$. For this let $\cB_3=\cB\sm\{ e_2\}$ and 
$\bG_3=R_{k/\Q}\bG_3'$ where $\bG_3'$ denotes the $k$-algebraic
group $SO(f_3)$ and $f_3$ is the restriction of the form $f$ to the subspace
of $E$ generated by $\cB_3$. It can be verified from the results of \cite{FOR}, section 2,
that there is an ideal $\ua_1$ of $\Z$ such that $P=X_3/\Phi_3$ is a connected closed
orientable totally geodesic $(n-1)$-dimensional submanifold of $X/\Phi$, where
$\Phi$ is any subgroup of finite index of $\Gamma (\ua_1 )$. Here $\Phi_3=\Phi\cap\bG_3(\Q )$
and $X_3=(K\cap G_3)/G_3$. Note that, since $\cB_2\sbs\cB_3$, we have that $T\sbs P$.\\

Now take an ideal $\ua '\sbs \ua_0 \cap \ua_1$. Define $\Phi '=\Phi (\ua ', \ub ', \ub ')$
where $\ub '=\ub (\ua ')$. Define $\Phi_1'$, $\Phi_2'$, $\Phi_3'$, $M'$, $N'$, $T'$, $P'$
accordingly. Then $M'$, $N'$, $T'$, $P'$ satisfy 1,2,3 of the statement of the Proposition.
That is, they satisfy all required properties, except the ones
about the widths of the normal geodesic tubular neighborhoods.
(Note that $N'$ and $P'$ intersect transversally because both are different connected 
totally geodesic hypersurfaces of $M'$.)\\

By using an argument similar to the one
in pp. 235-236 of \cite{FOR} we can find an ideal $\ua ''\sbs\ub '$ such that the following
holds. Define $\Phi ''=\Phi (\ua '', \ub '', \ub '')$
where $\ub ''=\ub (\ua '')$. Define $\Phi_1''$, $\Phi_2''$, $\Phi_3''$, $M''$, $N''$, $T''$, $P''$
accordingly. Then we can choose $\ua ''$ such that
the normal geodesic tubular neighborhood of $P''\cap N''$ is larger that $s$.
Note that $M''$, $N''$, $T''$, $P''$ also satisfy 1,2,3 of the statement of the Proposition.
(Here we required $\ua ''\sbs\ub '$ so that $\Phi''$ is a subgroup of $\Phi'$; follows that
$M''\ra M'$ is a finite cover.)\\

Finally, choose $\uc\sbs\ub ''$ as in \cite{FOR}, pp. 235-236,  such that the following
holds. Define $\Phi =\Phi (\ua '', \ub '', \uc)$ and
define $\Phi_1$, $\Phi_2$, $\Phi_3$, $M$, $N$, $T$, $P$
accordingly. Then we can choose $\uc$ such that
the normal geodesic tubular neighborhood of $N$ is larger that $r$. 
Note that we also have 
that $M$, $N$, $T$, $P$ satisfy 1,2,3 of the statement of the Proposition
and that the normal geodesic tubular neighborhood of $P\cap N$ is larger that $s$. 
(Note also that $N=N''$.) Then $M$, $N$, $T$, $P$ satisfy 1,2,3, 4 of the statement of the Proposition.
This completes the proof of the Proposition.
\vspace{.3in}

\noindent {\bf Remark.} To prove the assertion in the remark after the statement of the Proposition
just replace the roles of $N$ and $T$ in the proof above.\\

We now recall the construction of branched covers and introduce some notation
(for more details see \cite{Pan}).
Let $M$ be a hyperbolic manifold of dimension $n$ and let $R$ and $Q$ be compact
totally geodesic submanifolds of $M$ of dimensions $n-1$ and $n-2$ with
$\p Q=R$. Assume that the normal bundle of $Q$ is trivial. The $i$-branched cover $M(i)$
of $M$ with respect to $(Q,R)$ is obtained in the following way.
Cut $M$ along $Q$ to produce a manifold $M'$. Since the normal bundle of $Q$
is trivial $M'$ contains two copies of $Q$ that intersect in $R$. Label these
copies $Q_0$ and $Q_1$. Take now $i$ copies of $M'$ and identify $Q_1$ of the first copy of 
$M'$ with $Q_0$ of the second copy of $M'$, and so on. At the end 
identify $Q_1$ of the last copy of  $M'$ with $Q_0$ of the first copy of $M'$.
The resulting manifold is the $i$-branched cover of $M$ with respect to $(Q,R)$.
The branched cover $M(i)$ contains $i$ copies of $Q$. The union of these copies
is called the $i$-book of $M(i)$. These copies intersect in a (unique) copy
of $R$. Hence we consider $R$ as a submanifold of $M(i)$. There is a projection
map $\pi : M(i)\ra M$, which restricted to each copy of $M'$ in $M(i)$ is just the
identification map (identify back the two copies of $Q$). Note that if $R=\emptyset$
then $\pi :M(i)\ra M$ is an $i$-sheeted (ordinary) covering space. 
In general $\pi :M(i)\sm R\ra M\sm R$ is also an  $i$-sheeted covering space.\\

Gromov and Thurston proved that $M(i)$ admits a Riemannian metric with sectional curvatures
equal to -1 outside a normal geodesic tubular neighborhood of $R$  of width $s$, and sectional
curvatures in the interval $[-1-\epsilon ,-1]$, where $0< \epsilon <\frac{(10\, log\, 2i)^2}{s^2}$,
inside a normal geodesic tubular neighborhood of $R$. They also proved that, for large $i$ 
(how large depending on $M$ and $R$), $M(i)$ does not have that homotopy type of a locally
symmetric space (they mention that probably all $M(i)$, $i>1$ do not have that homotopy type of a locally
symmetric space). For the proofs of the Theorems
we will need the following stronger (and more technical) result, proved also by Gromov and Thurston
in section 3.5 of \cite{GT} (p. 6):\\

\noindent {\bf 2.1.} {\it  
Let $M_{s}$, $s =1,2,...$ be a sequence of closed $n$-dimensional hyperbolic manifolds,
$n>4$, and $R_{s}$ a closed codimension 2 totally geodesic submanifold of $M_{s}$ 
with trivial normal bundle that bounds
a compact codimension 1 totally geodesic submanifold $Q_{s}$ of $M_{s}$. Assume that 
the width of the normal geodesic tubular neighborhood of $R_{s}$ is larger than $s$.

Then there is a sequence of branched covers $M_{s}(i_{s})$ of $M_{s}$,
with respect to $(Q_{s},R_{s})$, satisfying the following property:

Given $\epsilon >0$ there is a $s_0$ such that for all $s\geq s_0$ we have: }
\begin{enumerate}
\item[{1.}] {\it $M_{s}(i_{s})$ does not have the homotopy type of a locally symmetric space.}

\item[{2.}] {\it $M_{s}(i_{s})$ admits a metric with sectional curvatures in $[-1-\epsilon , -1]$.}
\end{enumerate}

{\bf Remark.} Let $M$ be a closed negatively curved manifold of dimension $\neq$ 3,4.
Farrell-Jones Rigidity Theorem \cite{FJ2} implies that $M$ has the homotopy type of a locally
symmetric space if and only if $M$ admits a negatively curved locally symmetric differentiable structure,
i.e. a differentiable structure that supports a negatively curved locally symmetric Riemannian metric.
\vspace{.2in}

\noindent {\bf Proof of Theorem 1.} Fix $\epsilon >0$.
For each $s=1,2,3,...$ let 
$M_{s}'$, $N_{s}'$, $P_{s}'$ and $T_{s}'$ be as in the Proposition, 
with the width of the normal geodesic tubular neighborhood of $N_{s}'\cap P_{s}'$ larger
than $s$ and with large $r$
(how large will be determined in a moment; note that, by item 2 of the Proposition, $N_{s}'$ does not 
depend or $r$). We assume also $r>2s$. To alleviate the notation we will drop the subindex ``$s$'' whenever
this causes no confusion. Write $U'=N'\cap T'$.
By item 3 of the Proposition  $c'=Dual(U')\neq 0\in H^3(M',\Z_2 )$. Then we have a 
smooth structure $\Sigma '$ on $M'$ such that its corresponding $PL$ structure corresponds to
$c'$. (We choose this correspondence to assign the hyperbolic differentiable structure to $0\in H^3(M',\Z_2)$, 
see \cite{O}. It follows that $\Sigma '$ is not $PL$-concordant to the hyperbolic differentiable structure, see \cite{KS}.) 
Choose $r$ large enough so that $(M',\Sigma ')$ admits a Riemannian metric with sectional curvatures
$\epsilon$-close to -1, see \cite{O} ($r$ depends on $s$ and $\epsilon$). The important point here is 
the following. Identify (a piece of) the tubular neighborhood of $N'$ with $N'\times [-r,r]$.
Then outside $V'=N'\times [-r,-\frac{r}{2}]$ we have that $\Sigma '$ coincides with the hyperbolic 
differentiable structure and the $\epsilon$-pinched metric, mentioned above, is hyperbolic, see \cite{O}.
That is, the change of the differentiable structure and the change of the metric
only happen inside $V'$. Note that the distance between $V'$ and $N'$ is $\frac{r}{2}>s$.\\

\noindent {\bf Remark.} Since the differentiable structure $\Sigma '$ is not $PL$-concordant to 
the hyperbolic differentiable structure, we have that the identity map $M'\ra (M',\Sigma ')$
is not homotopic to a $PL$ homeomorphism. 
More generally, if $M$ is any closed nonpositively curved manifold of dimension $\neq$ 3,4, and
$\Sigma$ is a differentiable structure not $PL$-concordant to the given nonpositively
curved differentiable structure then the identity map $M\ra (M,\Sigma )$
is not homotopic to a $PL$ homeomorphism.
To see this, suppose that $h:M\times [0,1]\ra M$
is a homotopy from $id_{M}$ to a $PL$ homeomorphism. Then, by Farrell-Jones Rigidity Theorem
\cite{FJ2}, the map $H'(x,t)=(h(x,t),t)$ is homotopic to a homeomorphism $H:M\times [0,1]\ra M\times [0,1]$
with $H_0=id_M$ and $H_1$ equal to the $PL$ homeomorphism above. It follows that
$\Sigma$ is $PL$-concordant to the given nonpositively curved differentiable structure.\\

From item 3 of the Proposition it follows that
$0\neq [N']\in H_{n-1}(M',\Z_2 )$. 
Let $M$ be the double cover of $M'$ with respect
to $N'$, that is, with respect to the kernel of the map $\pi_1 M'\ra \Z_2$ given by
$\alpha\mapsto Dual(N')[\alpha ]$. Let $p:M\ra M'$ be the covering projection. Note that
$M$ consists of two copies $A$, $B$ of the manifold obtained from $M'$ by 
cutting along $N'$. $A$ and $B$ intersect in two copies of $N'$. Denote one of these copies
by $N$.
Let $T=p^{-1}(T')$. Write $U=N\cap T$. Then $p|_{U}:U\ra U'$ is a homeomorphism
(in fact an isometry). Hence $0\neq [U]\in H_{n-3}(M)$. 
Hence, as before,  $c=Dual(U)\neq 0\in H^3(M,\Z_2 )$.
Also, the normal geodesic tubular neighborhood of $N$ is as large as the normal geodesic
tubular neighborhood of $N'$ and let  $V$ be the piece of $p^{-1}(V')$ that lies
in the $r$ normal geodesic tubular neighborhood of $N$. Also we assume that, say, $V\sbs A$.
Note that $V$ is a piece, of width $\frac{r}{2}$, of the normal geodesic tubular neighborhood of $N$
which is at a distance $\frac{r}{2}>s$ from $B$.\\

Let $P=p^{-1}(P')$, $Q=P\cap B$ and $R=\p Q$. Note that
$R$ consists of two copies of $N'\cap P'$ and that the width of the normal
geodesic tubular neighborhood of $R$ is larger that $s$. 
Note also that the distance between $V$ and $Q$ is $>s$.\\

So far we have 
obtained (we write back the subindex ``$s$'') sequences $M_{s}$, $N_{s}$, $T_{s}$,
$R_{s}$, $Q_{s}$, $V_{s}$, $U_{s}=N_s\cap T_s$, $c_{s}$ satisfying: 
\begin{enumerate}
\item[{a.}]  $M_{s}$ is a closed $n$-dimensional hyperbolic manifold
and $N_{s}$ is a closed codimension one totally geodesic submanifold
of $M_{s}$. $U_{s}$ is a closed codimension three totally geodesic submanifold
of $M_{s}$ and $U_{s}\sbs N_{s}$.

\item[{b.}] $R_{s}$ is a closed codimension 2 totally geodesic submanifold of $M_{s}$ that
bounds the compact codimension 1 totally geodesic submanifold $Q_{s}$. 
The width of the normal geodesic tubular neighborhood of $R_{s}$ is larger than $s$.

\item[{c.}] $V_{s}$ is a piece of the tubular neighborhood of $N_{s}$ and $V_{s}$ is at a distance
$>s$ from $Q_{s}$.

\item[{d.}] $c_{s}=Dual(U_{s})$ is a nonzero cohomology class. As for $M_s'$, $M_{s}$, with the smooth structure
$\Sigma _s$ whose corresponding $PL$ structure corresponds to $c_{s}$, admits a Riemannian metric with sectional curvatures
$\epsilon$-close to -1. Outside $V_{s}$  this $\epsilon$-pinched to -1 metric is hyperbolic.
\end{enumerate}

By 2.1, there is a $s_0$, such that for all $s\geq s_0$ we have: 
\begin{enumerate}
\item[{1.}]  $M_{s}(i_{s})$ does not have the homotopy type of a locally symmetric space.

\item[{2.}] $M_{s}(i_{s})$ admits a metric with sectional curvatures in $[-1-\epsilon , -1]$.
This metric is hyperbolic outside the normal geodesic tubular neighborhood of width $s$ of
$R_s$.
\end{enumerate}

Write $(M_1)_{s}=M_{s}(i_{s})$ and we again drop the subindex ``$s$''. Let $\pi :M_1\ra M$
be the ramified projection. Choose a copy of $M\sm ( Q\sm R)$ in $M_{1}$. Hence we then can find
$N_{1}$, $V_{1}$, $U_{1}$ contained in this copy such that:
\begin{enumerate}
\item[{}] $\pi|_{N_{1}}:N_1\ra N$, $\pi|_{V_{1}}:V_1\ra V$, $\pi|_{U_{1}}:U_1\ra U$
are homeomorphisms. 
\end{enumerate}

Since  $U$ does not bound, it follows that $c_1=Dual(U_1)\neq 0\in H^3(M_1,\Z_2)$.
Then we have a 
smooth structure $\Sigma_1$ on $M_1$ such that its corresponding $PL$ structure corresponds to
$c_1$. (Again, we choose this correspondence to assign the branched differentiable structure to $0\in H^3(M',\Z_2)$.)
Hence $\Sigma_1$ is not $PL$-concordant to the branched differentiable structure.
But $M_2=(M_1,\Sigma_1)$ admits a Riemannian metric with sectional curvatures
$\epsilon$-close to -1 such that outside $V_1$  this metric coincides with the metric of
item 2 above. It follows then that the identity map $id :M_1\ra M_2$ is a homotopy equivalence
between $\epsilon$-pinched to -1 closed manifolds which is not homotopic
to a $PL$ homeomorphism (see the last remark above)
and $M_1$ and $M_2$ do not have the homotopy type of a locally symmetric space.
This proves Theorem 1.
\vspace{.3in}

Theorem 2 can be directly deduced from Theorem 1 and the $C^\infty -$ Hauptvermutung of Scharlemann
and Siebenmann \cite{SC1}. For more details see \cite{FJO2} or \cite{FO3}.\\

Before we prove Theorem 3, we need a remark about coverings of branched covers. Let $p:M\ra M'$ be a
cover, where $M$ and $M'$ are hyperbolic manifolds. Let $R$ and $R'$ be closed codimension two
totally geodesic submanifolds of $M$ and $M'$, respectively, that bound
closed codimension one totally geodesic submanifolds $Q$ and $Q'$ of $M$ and $M'$, respectively. Assume $Q=p^{-1}(Q')$
and $R=p^{-1}(R')$.
Then, for any $i$ we can use $p$ to construct a cover $q: M(i)\ra M'(i)$, where $M(i)$ and $M'(i)$ are
the $i$-branched covers of $M$ and $M'$ with respect to $(Q,R)$ and $(Q',R')$ respectively.
The covers $p$ and $q$ fit in the following commutative square:

$$
\begin{array}{ccc} M(i) & \stackrel{\pi}{\rightarrow} & M \\
q\downarrow & &\downarrow $p$\\
 M'(i) & \stackrel{\pi '}{\rightarrow} & M' 
 \end{array}
$$

If $p$ is a $\ell$-sheeted cover, then so is $q$. 
If $p$ is a regular cover, so is $q$.
We now prove Theorem 3.\\

\noindent {\bf Proof of Theorem 3.}
Let $M'$, $N'$, $T'$, $R'$, $Q'$, $V'$, $U'$, $c'$, $\Sigma '$
be as in the proof of Theorem 1 (we are dropping the subindex ``$s$'' and introducing
a prime on each symbol). These objects satisfy properties a,b,c,d in the proof of
Theorem 1. We assume also $T'\sbs P'$
(see remark after the statement of the Proposition).
Let also $M_1'$, $N_{1}'$, $V_{1}'$, $U_{1}'$, $c_1'$, $\pi '$, $\Sigma_1 '$, $M_2'$ 
be as in the proof of Theorem 1 (we are also dropping the subindex ``$s$'' and introducing
a prime on each symbol).
We have that $M_1'$ satisfies properties a and b in the proof of Theorem 1 (assuming $s$ large enough).\\

Let $p:M\ra M' $ be the double cover of $M'$ with respect to $N'$. 
Note that $M$ consists of two copies $D$, $E$ of the manifold obtained from $M'$ by 
cutting along $N'$. $D$ and $E$ intersect in two copies of $N'$.
Let $\Sigma =p^*\Sigma '$.
We know from \cite{FOR} that the identity $M\ra (M,\Sigma )$ is now
homotopic to a diffeomorphism $f:M\ra (M,\Sigma )$ but the unique harmonic map $h$
homotopic to $f$ (or to the identity) is not one-to-one. This is because the $PL$
structure corresponding to $\Sigma$ corresponds to $c=p^*(c')$ which vanishes.
(The cohomology class $c$ vanishes because it is dual to $p^{-1}(U')$, which
is the boundary of $D\cap p^{-1}(T')$.)\\

Let $M_1=M(i)$. Since we also have (by definition) $M_1'=M'(i)$ the commutative square
given just before this proof becomes:

$$
\begin{array}{ccc} M_1 & \stackrel{\pi}{\rightarrow} & M \\
q\downarrow & &\downarrow $p$\\
 M_1' & \stackrel{\pi '}{\rightarrow} & M' 
 \end{array}
$$

Define $\Sigma_1=q^*\Sigma_1'$, $M_2=(M_1,\Sigma_1)$, $c_1=q^*c_1'$ and $U_1=q^{-1}(U_1')$.
Note that:

\begin{enumerate}
\item[{(i)}] $q: M_2\ra M_2'$ is a smooth cover.

\item[{(ii)}] $M_1'=M_1/F$, where 
$F=\{ id, \phi\}\cong\Z_2$ and $\phi :M_1\ra M_1$ is the unique nontrivial covering
transformation.

\item[{(iii)}] The $PL$ structure corresponding to $\Sigma_1$
corresponds to $c_1$.

\item[{(iv)}] $c_1$ is dual to $U_1$.
\end{enumerate}

\noindent {\it Claim.} $M_1$ is not homotopy equivalent to a locally symmetric space.
\vspace{.07in}

\noindent If $M_1$ supports a negatively curved locally symmetric differentiable structure (see remark before
the proof of Theorem 1) then, by Mostow's Rigidity Theorem, $\phi$ can be realized by 
an isometry. It follows that $M_1'$ admits a negatively curved locally symmetric differentiable structure.
This is a contradiction. This proves the claim.\\

Now, since $p: M_1\ra M_1'$  and $q: M_2\ra M_2'$ are smooth covers we have that
both $M_1$, $M_2$ admit $\epsilon$-pinched to -1 Riemannian metrics.\\

It remains to prove that $\Sigma_1$ is $DIFF$-concordant to the differentiable structure
of $M_1$. To be able to repeat the argument given in \cite{FOR}, pp. 230-233, we need to prove
that the cohomology class $c_1$ vanishes. Equivalently we need to prove that $U_1$ bounds. 
Recall that $U'=T'\cap N'$. We assume that $T'\sbs  P'$ (see remark
after the statement of the Proposition).
Then we can also find $T_1'$ in one of the copies that form $M_1'$ such that
$\pi '|_{T_1'}:T_1'\ra T'$ is a homeomorphism and $U_1'\sbs T_1'$. Then $U_1'=T_1'\cap N_1'$.
Let $T_1=q^{-1}(T_1')$. A simple geometric argument then shows that $U_1=\p (T_1\cap D_1)$,
where $D_1$ lies in the same copy where $T_1$ lies and $\pi (D_1)=D$. This proves 
Theorem 3.
\vspace{.3in}

\noindent {\bf Proof of Theorem 4.} We use all notation from the proof of Theorem 3, with the following
changes. 
\begin{enumerate}
\item[{1.}] Now we assume that the width of the normal geodesic tubular neighborhood of $T'$ (not of $N'$)
is larger than $r$ (see remark after the statement of the Proposition). It follows that the 
width of the normal geodesic tubular neighborhood of $T$ is also larger than $r$.

\item[{2.}] The changes of structure and metric happen now in a piece of the normal geodesic
tubular neighborhood of $T'$ and $T$: since the normal bundle of $T'$ is trivial we have
${\cal{N}}_r(T')\sm T'$ can be identified with $T'\times \bS^1\times (0,r)$, where  
${\cal{N}}_r(T')$ is the  normal geodesic tubular neighborhood of $T'$ of width $r$.
In \cite{FO1} it is shown that we can take now $V'=T'\times I\times (0,\frac{r}{2})$, where
$I\subset\bS^1$ is any non-trivial interval. That is, outside $V'$ the differentiable structure $\Sigma '$ 
coincides with the hyperbolic differentiable structure, and the metric is hyperbolic.
Note that, by choosing $I$ properly, we have that $V'$ does not intersect $P'$.

\item[{3.}] Now $p:M\ra M'$ denotes the finite sheeted cover given in \cite{FO1}. Again we have that
$\Sigma $ is now $DIFF$-concordant to the hyperbolic differentiable structure. The new feature now
is that the metric pulled back from $(M',\Sigma ')$ can be deformed to the hyperbolic one and all
this deformation happens inside $V$.
\end{enumerate}

Now, since $V$ does not intersect $P$, we can also deform the metric of $M_1$ to the metric pulled
back from $M_2$, and all this happens inside $V_1$. This is the main ingredient needed for 
the proof. The rest follows exactly as in \cite{FO1}. This completes the proof of Theorem 4.
\vspace{.3in}

{\bf Proof of Theorem 5.}  The proof follows the ideas of the proofs of Theorem and Corollary 3 in
\cite{FO2}. We use all objects in the proof of Theorem 4 above, but change the notation a bit
to match the notation in \cite{FO2}. From the proof above we have a homotopy commutative diagram:

$$
\begin{array}{ccc} M_1 & \stackrel{f}{\rightarrow} & M_2 \\
q\downarrow & &\downarrow p\\
 M_1' & \stackrel{f '}{\rightarrow} & M_2' 
 \end{array}
$$

\noindent where $f$ is a diffeomorphism and $f'$ is not homotopic to a $PL$ homeomorphism. 
Recall that $f'$ is the identity and $f$ is homotopic to the identity
(the underlying topological manifolds of $M_2'=(M_1',\Sigma ')$ and $M_2=(M_1,\Sigma )$
are $M_1'$ and $M_1$ respectively).
We now change the notation to match the one in \cite{FO2}.
First identify $M_1$ with $M_2$ via $f$ and write $N$. Then
write $M_0$ for $M_1'$, $M_1$ for $M_2'$, $f$ for $f'$, $p_0$ for $q$ and $p_1$ for
$pf$. The diagram above becomes:

$$
\begin{array}{ccc}  & N &  \\
p_0\swarrow & &\searrow p_1\\
 M_0 & \stackrel{f }{\rightarrow} & M_1 
 \end{array}
$$

As in \cite{FO2} we have then that $M_0$, $M_1$, $N$, are smooth closed manifolds
of dimension $n>10$, that admit Riemannian metrics $g_0$, $g_1$ on $M_0$
and $M_1$, respectively, and smooth regular finite covers 
$p_0 : N\rightarrow M_0$,   $p_1 : N\rightarrow M_1$ such that: \\

(1) The map $f$ is not homotopic to a $PL$-homeomorphism.

(2) and $g_0$ and $g_1$ have sectional curvatures in $[-1-\epsilon ,-1+\epsilon]$.

(4) There is a $C^\infty$ family of $C^\infty$ Riemannian metrics $h_s$ on $N$,
$0\leq s\leq 1$, with $h_0 =p_0 ^* g_0$, and $h_1 =p_1 ^* g_1$,
such that every $h_s$ has sectional curvatures in $[-1-\epsilon ,-1+\epsilon]$.\\

Let $G_i\subset Diff\, (N)$, be (finite) subgroups of the group $Diff\, (N)$, of all self-diffeomorphisms
of $N$, such that $N/G_i=M_i$, $i=0,1$. 
It was shown in \cite{FOR}, \cite{FO1}, that  $G_0$ and $G_1$ are conjugate in $Top\, N$, 
via a homeomorphism homotopic to $id_N$; hence $\gamma G_0=\gamma G_1$, where
$\gamma :Diff \, (N) \ra Out\, (\pi_1 N)$ is the map described in the proof of Corollary 3 in
\cite{FO2}. Note that $G_i\subset Iso(N, h_i)$,  
where $Iso\, (N, h_i)\subset Diff\, (N)$ is the subgroup consisting 
of all isometries of the negatively curved  manifold $(N, h_i)$.\\

If the Ricci flow does not converge smoothly for some $h_s$, we are done. So, let us assume
that the Ricci flow converges smoothly for all $h_s$. We will show a contradiction.
Write $h_{s,t}$, for the Ricci flow starting at $h_{s,0}=h_s$, $0\leq t <\infty$, converging to the
negatively curved Einstein metric $j_{s}$. Using the same argument as the one given
in the proof of Theorem of \cite{FO2} we get that all $j_{s}$ are equal (modulo
diffeomorphism and rescaling). Moreover, there is a diffeomorphism $\phi: N\ra N$ homotopic to the
identity such that $j_0$ is equal to $\phi^*j_1$ up to scaling (see also the proof of Corollary 3 in \cite{FO2}). 
It follows that  $G_1$ is conjugate in $Diff\, (N)$
to a subgroup of $Iso\, (N,j_0)$ 
via a diffeomorphism $\phi$ homotopic to $id_N$; i.e. $\phi^{-1}G_1\phi \subset Iso\, (N,j_0)$. 
Note that  $\gamma (\phi^{-1}G_1 \phi )=\gamma (G_1)$ since $\phi\sim id_N$; hence 
$\gamma (\phi^{-1}G_1 \phi )=\gamma (G_0)$. This implies that  $\phi^{-1}G_1 \phi =G_0$ since both
$\phi^{-1}G_1 \phi$ and $G_0$ are subgroups of $Iso\, (N,j_0)$ and Borel-Conner-Raymond
showed (see \cite{CR}, p.43)
that $\gamma$ restricted to compact subgroups of $Diff\, (N)$ is monic. (Recall that $N$ is aspherical
and the center of $\pi_1(N)$ is trivial.)
It follows that $\phi$ induces a diffeomorphism $\varphi$ between $M_0=N/G_0$ and 
$M_1=N/G_1$. To find a contradiction we have to prove that $\varphi$ is homotopic to the
identity. Since all manifolds here a aspherical, it is enough to prove that the induced
map $\varphi_* $ at the fundamental group level is the identity. Note that, since $\phi_*$ 
is the identity, we have that $\varphi_* |_{\pi_1 N}$ is the identity. Then the fact that
$\varphi_*$ is also the identity follows from the next Lemma:\\

{\bf Lemma.} {\it Let $\Gamma $ be the fundamental group of a closed negatively curved manifold.
 Let $H$ be a subgroup of finite index of $\Gamma$
and let $\alpha :\Gamma\ra \Gamma$ be an isomorphism whose restriction to $H$ is
the identity. Then $\alpha$ is also the identity.}\\

{\bf Proof.} 
Let $x\in \Gamma$. Since $H$ has finite index in $\Gamma$ we have that there is a positive
integer $n$ such that  $x^n\in H$. Hence $\alpha (x^n)=x^n$. Therefore $(\alpha (x))^n=x^n$.
But $n$-roots are unique in $\Gamma$, thus $\alpha (x)=x$. This proves the Lemma.\\

F.T. Farrell

SUNY, Binghamton, N.Y., 13902, U.S.A.\\

P. Ontaneda

UFPE, Recife, PE 50670-901, Brazil

\end{document}